\documentclass{article}

\usepackage{amssymb}
\usepackage{amsfonts}
\usepackage{graphicx}
\usepackage{epsfig}

\usepackage[all]{xy}

\newtheorem{prop}{Proposition}[section]
\newtheorem{teo}[prop]{Theorem}
\newtheorem{remark}[prop]{Remark}
\newtheorem{lemma}[prop]{Lemma}

\begin{document}

\setlength{\parindent}{0mm}

\author{\small{\textbf{ Marian Ioan MUNTEANU }}}
\title{\textbf{Cheeger Gromoll type metrics on the tangent bundle}}
\date{}
\maketitle

\def\proof{{\sc Proof.\ }}
\newcommand{\gata}{\hfill\hskip -1cm \rule{.5em}{.5em}}
\def\stackunder#1#2{\mathrel{\mathop{#2}\limits_{#1}}}

\newcommand\uu{\mathtt{u}}
\def\RR{{\mathcal{R}}}
\def\LL{\mathtt{L}}
\def\MM{\mathtt{M}}
\def\NN{\mathtt{N}}
\def\kk{\mathtt{k}}
\def\rr{{\mathfrak r}}

\begin{abstract}
\noindent
In this paper we study a Riemanian metric on the tangent bundle
$T(M)$ of a Riemannian manifold $M$ which generalizes the Cheeger Gromoll
metric and a compatible almost complex structure
which together with the metric confers to $T(M)$ a structure of
locally conformal almost K\"ahlerian manifold.
We found conditions under which $T(M)$ is almost
K\"ahlerian, locally conformal K\"ahlerian or K\"ahlerian or when
$T(M)$ has constant sectional curvature or constant scalar
curvature.
\end{abstract}

\section[]{Preliminaries}

Given a Riemannian manifold $(M,g)$ one can define
several Riemannian metrics on the tangent bundle $T(M)$ of $M$.
Maybe the best known example is the Sasaki metric $g_S$ introduced
in \cite{kn:Sas58}. Although the Sasaki metric is {\em naturally}
defined, it is very {\em rigid}. For example, the Sasaki metric is
not, generally, Einstein. Or, the tangent bundle $T(M)$
with the Sasaki metric is never locally symmetric unless the metric
$g$ on the base manifold is flat (see  \cite{kn:Kow71}). E.Musso \& F.Tricerri
\cite{kn:MT88} have proved that the Sasaki
metric has constant scalar curvature if and only if $(M,g)$ is
locally Euclidian. In the same paper, they have given an explicit
expression of a positive definite Riemannian metric introduced by
J.Cheeger and D.Gromoll in \cite{kn:CG72} and called this metric
{\em the Cheeger-Gromoll metric}. M.Sekizawa (see \cite{kn:Sek91}),
computed geometric objects related to this metric.
Later, S.Gudmundson and E.Kappos in
\cite{kn:GK02}, have completed these results and have shown that the scalar curvature of
the Cheeger Gromoll metric is never constant if the metric on the
base manifold has constant sectional curvature. Furthermore,
M.T.K.Abbassi \& M.Sarih have proved that $T(M)$ with the Cheeger
Gromoll metric is never a space  of constant sectional curvature
(cf. \cite{kn:AS03}). It is also known that the tangent bundle $T(M)$
of a Riemannian manifold $(M,g)$ can be organized as an almost K\"ahlerian
manifold (see \cite{kn:Dom62}) by using the decomposition of the tangent
bundle to $T(M)$ into the vertical and horizontal distributions, $VTM$
and $HTM$ respectively (the last one being defined by the Levi Civita
connection on $M$), the Sasaki metric and an almost complex structure
defined by the above splitting. A more general metric is given by M.Anastasiei
in \cite{kn:Ana99} which generalizes both of the Sasaki and Cheeger Gromoll
metrics: it preserves the
orthogonality of the vertical and horizontal distributions, on the horizontal
distribution it is the same as on the base manifold, and finally the
Sasaki and the Cheeger Gromoll metric can be obtained as particular
cases of this metric. A compatible almost complex structure is also
introduced and hence $T(M)$ becomes a locally conformal almost
K\"aherian manifold. On the other hand,
V.Oproiu and his collaborators constructed a family of Riemannian
metrics on the tangent bundles of Riemannian manifolds which possess
interesting geometric properties (cf. \cite{kn:Opr99, kn:Opr199,
kn:Opr01, kn:OP04}) (for example, the scalar curvature of $T(M)$
can be constant also for a non-flat base manifold with constant
sectional curvature). Then M.T.K.Abbassi \& M.Sarih proved in
\cite{kn:AS05} that the considered metrics by Oproiu form a
particular subclass of the so-called $g$-{\em natural metrics} on
the tangent bundle (see also \cite{kn:Aba04, kn:AS05,
kn:AS105, kn:AS205, kn:KS88}).\\[3mm]
In this paper we described a family $g_a$ of Riemannian metrics
of Cheeger Gromoll type, on the tangent bundle $T(M)$ of the Riemannian manifold
$(M,g)$ and a compatible almost complex structure $J_a$ which
bestow to $T(M)$ a structure of locally conformal almost K\"ahlerian manifold.
We found an almost K\"ahlerian structure on $T(M)$ and we proved that there is no
Cheeger Gromoll type structure on $T(M)$ such that the manifold $(T(M),g_a,J_a)$
is K\"ahlerian. We studied the possibility for the sectional curvature on $T(M)$ to
be constant and we found a flat metric on $T(M)$ (of course of Cheeger Gromoll type).
Finally, if $M$ is a real space form, we were interested to find when $T(M)$
endowed with the metric $g_a$ has constant scalar curvature.

\section{On the Geometry of the Tangent Bundle $T(M)$}

Let $(M,g)$ be a Riemannian manifold and let $\nabla$ be its Levi
Civita connection. Let $\tau:T(M)\longrightarrow M$ be the tangent
bundle. If $\uu\in T(M)$ it is well known the following decomposition of the
tangent space $T_\uu T(M)$ (in $\uu$ at $T(M)$)
$$
    T_\uu T(M)=V_\uu T(M)\oplus H_\uu T(M)
$$
where $V_\uu T(M)=\ker \tau_{*,\uu}$ is the vertical space and
$H_\uu T(M)$ is the horizontal space in $\uu$ obtained by using $\nabla$.
(A curve $\widetilde\gamma:I\longrightarrow T(M)$ , $t\mapsto(\gamma(t),V(t))$
is {\em horizontal} if the vector field $V(t)$ is parallel along
$\gamma=\widetilde\gamma\circ\tau$. A vector on $T(M)$ is {\em horizontal}
if it is tangent to an horizontal curve and {\em vertical} if it is tangent
to a fiber. Locally, if $(U,x^i)$, $i=1, \ldots, m$, $m=\dim M$,
is a local chart in $p\in M$, consider
$(\tau^{-1}(U), x^i, y^i)$ a local chart on $T(M)$. If $\Gamma_{ij}^k(x)$ are the
Christoffel symbols, then
$\delta_i=\frac\partial{\partial x^i}-\Gamma_{ij}^k(x)y^j\ \frac\partial{\partial y^k}$
in $\uu$, $i=1,\ldots,m$ span the space $H_\uu T(M)$, while
$\frac\partial{\partial y^i}$, $i=1,\ldots,m$ span the vertical space $V_\uu T(M)$.)
We have obtained the horizontal (vertical) distribution $HTM$ ($VTM$) and a direct sum decomposition
$$
   TTM=HTM\oplus VTM
$$
of the tangent bundle of $T(M)$.
If $X\in\chi(M)$, denote by $X^H$ (and $X^V$, respectively) the horizontal lift
(and the vertical lift, respectively) of $X$ to $T(M)$.

If $\uu\in T(M)$ then we consider the energy density in $\uu$ on $T(M)$,
namely $t=\frac12\ g_{\tau(\uu)}(\uu,\uu)$.

\subsection{The Cheeger-Gromoll Structure}

The Cheeger-Gromoll metric on $T(M)$ is given by
\begin{equation}
\left\{
\begin{array}{l}
{g_{CG}}_{(p,\uu)}(X^H,Y^H)=g_p(X,Y),\
{g_{CG}}_{(p,\uu)}(X^H,Y^V)=0\\
\displaystyle
{g_{CG}}_{(p,\uu)}(X^V,Y^V)=\frac1{1+2t}\ \left(g_p(X,Y)+g_p(X,\uu)g_p(Y,\uu)\right)
\end{array}
\right.
\end{equation}
for any vectors $X$ and $Y$ tangent to $M$.
Moreover, an almost complex structure $J_{CG}$, compatible
with the Chegeer-Gromoll metric, can be defined by the formulas
\begin{equation}
\left\{
\begin{array}{l}
J_{CG}X^H_{(p,\uu)}=\rr X^V-\frac1{1+\rr}\ g_p(X,\uu)\uu^V\\[2mm]
J_{CG}X^V_{(p,\uu)}=-\frac1\rr\ X^H-\frac1{\rr(1+\rr)}\ g_p(X,\uu)\uu^H
\end{array}
\right.
\end{equation}
where $\rr=\sqrt{1+2t}$ and $X\in T_p(M)$. Remark that $J_{CG}\uu^H=\uu^V$
and $J_{CG}\uu^V=-\uu^H$. We get an almost Hermitian manifold $(T(M),J_{CG},g_{CG})$.
If we denote by $\Omega_{CG}$ the Kaehler 2-form (namely
$\Omega_{CG}(U,V)=g_{CG}(U,J_{CG}V),$ $\forall U,V\in\chi(T(M)))$ one
can prove the following
\begin{prop}
We have
\begin{equation}
d\Omega_{CG}=\omega\wedge\Omega_{CG},
\end{equation}
where $\omega\in\Lambda^1(T(M))$ is defined by
$$\omega_{(p,\uu)}(X^H)=0 {\ \rm and\ }
  \omega_{(p,\uu)}(X^V)=-\left(\frac1{\rr^2}+\frac1{1+\rr}\right)g_p(X,\uu),
 X\in T_p(M).
$$
\end{prop}
\proof
A simple computation gives the differential of $\Omega_{CG}$:

\quad $d\Omega_{CG}(X^H,Y^H,Z^H)=d\Omega_{CG}(X^H,Y^H,Z^V)=d\Omega_{CG}(X^V,Y^V,Z^V)=0$

\quad $d\Omega_{CG}(X^H,Y^V,Z^V)=\frac1\rr\left(\frac1{\rr^2}+\frac1{1+\rr}\right)
        \left[g(X,Y)g(Z,\uu)-g(X,Z)g(Y,\uu)\right]$

for any $X,Y,Z\in\chi(M)$.
Hence the statement. \gata

\begin{remark}\rm
The almost Hermitian manifold $(T(M),J_{CG},g_{CG})$ is never almost Kaehlerian (i.e. $d\Omega_{CG}\neq0$).
\end{remark}

Finally, a necessary condition for the integrability of $J_{CG}$  is that the base manifold
$(M,g)$ is locally Euclidian.

\vspace{2mm}

\subsection{The Cheeger Gromoll Type Structure}

A general metric, let's call it $g_a$, is in fact a family of Riemannian
metrics, depending on a parameter $a$, and the Cheeger-Gromoll
metric is obtained by taking $a(t)=\frac 1{1+2t}$. It is defined by the following formulas
(see also \cite{kn:Ana99})
\begin{equation}
\left\{
\begin{array}{l}
{g_a}_{(p,\uu)}(X^H,Y^H)=g_p(X,Y)\\
{g_a}_{(p,\uu)}(X^H,Y^V)=0\\
{g_a}_{(p,\uu)}(X^V,Y^V)=a(t)\big( g_p(X,Y)+g_p(X,\uu)g_p(Y,\uu)\big),
\end{array}
\right.
\end{equation}
for all $X,Y\in\chi(M)$, where $a:[0,+\infty)\longrightarrow(0,+\infty)$.

\begin{prop} {\rm (see also \cite{kn:Mun06})}
The metric defined above can be construct by using the method described by Musso
and Tricerri in {\rm \cite{kn:MT88}}.
\end{prop}

We intend to find an almost complex structure on $T(M)$, call it $J_a$, compatible
with the metric $g_a$. Inspired from the previous cases we look for the almost complex structure $J_a$
in the following way
\begin{equation}
\left\{
\begin{array}{l}
J_a X^H_{(p,\uu)}=\alpha X^V + \beta g_p(X,\uu)\uu^V\\[2mm]
J_a X^V_{(p,\uu)}= \gamma X^H + \rho g_p(X,\uu)\uu^H
\end{array}
\right.
\end{equation}
where $X\in\chi(M)$ and $\alpha$, $\beta$, $\gamma$ and $\rho$ are smooth functions on $T(M)$
which will be determined from $J_a^2=-I$ and from the compatibility conditions with the metric
$g_a$. Following the computations made in \cite{kn:Ana99} we get first $\alpha=\pm\frac1{\sqrt a}$
and $\gamma=\mp \sqrt{a}$.\\[2mm]
Without lost of the generality we can take
$\alpha=\frac1{\sqrt a} {\quad\rm  and\quad} \gamma=- \sqrt{a}.$
Then one obtains
$
\begin{array}{l}
\beta=-\frac1{2t}\left(\frac1{\sqrt{a}}+\epsilon\frac1{\sqrt{a+2bt}}\right) {\ \rm and\ }
\rho=\ \frac1{2t}\left(\sqrt{a}+\epsilon\sqrt{a+2bt}\right)
\end{array}
$
where $\epsilon=\pm1$.
\begin{remark}\rm
In this general case $J_a$ is defined on $T(M)\setminus 0$ (the bundle of non zero tangent vectors),
but if we consider $\epsilon=-1$ the previous relations define $J_a$ on all $T(M)$.
From now on we will work with $\epsilon=-1$.
\end{remark}
We have the almost complex structure $J_a$
\begin{equation}
\left\{
\begin{array}{l}
J_aX^H=\frac1{\sqrt{a}}\left(X^V-\frac1{\rr(1+\rr)}\ g(X,\uu)\uu^V\right)\\[3mm]
J_aX^V=-\sqrt{a}\left(X^H+\frac1{1+\rr}\ g(X,\uu)\uu^H\right).
\end{array}
\right.
\end{equation}
One obtains an almost Hermitian manifold $(T(M),g_a,J_a)$.\\[2mm]
If we denote by $\Omega_a$ the K\"ahler 2-form, $\Omega_a(U,V)=g_a(U,J_aV),$ $\forall U,V\in\chi(T(M))$
one obtains
\begin{prop}
{\rm (see also \cite{kn:Ana99})}
The almost Hermitian manifold $(T(M),g_a,J_a)$ is locally conformal almost K\"ahlerian, that is
\begin{equation}
d\Omega_a=\omega\wedge\Omega_a
\end{equation}
where $\omega$ is a closed and globally defined $1-$form on $T(M)$ given by
$$\omega(X^H)=0 {\quad and \quad}
\omega(X^V)=\left(\frac{a'}{a}-\frac1{1+\rr}\right)g(X,\uu).$$
\end{prop}

As consequence one can state the following

\begin{teo}
\label{th:aK}
The almost Hermitian manifold $(T(M),g_a,J_a)$ is almost K\"ahlerian if and only if
\begin{equation}
\label{eq:a}
\begin{array}{l}
a(t)=const\cdot\frac{e^{\sqrt{1+2t}}}{1+\sqrt{1+2t}}\ .
\end{array}
\end{equation}
\end{teo}
\proof
The result is obtained by integrating the equation
$\frac{a'}{a}=\frac1{1+\rr}$.
\gata

We will take $a(\rr)=\frac{2e^{\rr-1}}{1+\rr}$ if we ask $a(0)=1$.

\subsection{The Integrability of $J_a$.}

In order to have an integrable structure $J_a$ on $T(M)$ we have to compute the Nijenhuis
tensor $N_{J_a}$ of $J_a$ and to ask that it vanishes identically.\\[2mm]
For the integrability tensor $N_{J_a}$ we have the following
relations\\[1mm]
$
N_{J_a}(X^H,Y^H)=\displaystyle\frac{2a-(1+\rr)a'}{2a^2\rr(1+\rr)}\ \big(g(X,\uu)Y-g(Y,\uu)X\big)^V+(R_{XY}\uu)^V
$\\[1mm]
$
N_{J_a}(X^V,Y^V)=\big(-aR_{XY}\uu-\frac{a}{1+\rr}\ g(Y,\uu)R_{X\uu}\uu+\frac{a}{1+\rr}\ g(X,\uu)R_{Y\uu}\uu\big)^V-
$
\begin{equation}
\label{eq:NIJ}\begin{array}{l}
 -\left(\frac {a'}{2a}-\frac1{1+\rr}\right)\big(g(Y,\uu)X-g(X,\uu)Y\big)^V.
 \end{array}
\end{equation}
The expression for $N_{J_a}(X^H,Y^V)$ is very complicated.

Thus if $J_a$ is integrable then
$$
R_{XY}\uu=\frac{2a-(1+\rr)a'}{2a^2\rr(1+\rr)}\ \big(g(Y,\uu)X-g(X,\uu)Y\big)
$$
for every $X,Y\in\chi(M)$ and for every point $\uu\in T(M)$. It follows that $M$ is a space form $M(c)$
($c$ is the constant sectional curvature of $M$).
Consequently,
$$
a(\rr)=\frac{e^{2\rr}}{(1+\rr)\left(ce^{2\rr}(\rr-1)+k(1+\rr)\right)}
$$
with $k$ a positive real constant and $c$ must be nonnegative.

\vspace{2mm}

{\bf Question:} {\em Can $(T(M),g_a,J_a)$ be a K\"ahler manifold?}\\[2mm]
If this happens then we have to find an appropriate constant in
(\ref{eq:a}) such that the expression
$\frac{2a-(1+\rr)a'}{2a^2\rr(1+\rr)}$ is also a constant.

\begin{teo}
\label{th:K}
There is no Cheeger Gromoll type structure on $T(M)$ such that the manifold $(T(M),g_a,J_a)$
is K\"ahlerian.
\end{teo}
Now we give
\begin{prop}
Let $(M,g)$ be a Riemannian manifold and let $T(M)$ be its tangent bundle equipped with the metric $g_a$.
Then, the corresponding Levi Civita connection $\tilde\nabla^a$ satisfies the following relations:
\begin{equation}
\left\{
\begin{array}{l}
\tilde\nabla^a_{X^H}Y^H=(\nabla_XY)^H-\frac12\ (R_{XY}\uu)^V\\[2mm]
\tilde\nabla^a_{X^H}Y^V=(\nabla_XY)^V+\frac a2\ (R_{\uu Y}X)^H\\[2mm]
\tilde\nabla^a_{X^V}Y^H=\frac a2\ (R_{\uu X}Y)^H\\[2mm]
\tilde\nabla^a_{X^V}Y^V=\LL \left(g(X,\uu)Y^V+g(Y,\uu)X^V\right)+\frac{1-L}{\rr^2}\ g(X,Y)\uu^V-\\[1mm]
   \qquad\qquad\qquad\qquad      -\frac{L}{\rr^2}\ g(X,\uu)g(Y,\uu)\uu^V,
\end{array}
\right.
\end{equation}
where $\LL=\frac{a'(t)}{2a(t)}$.
\end{prop}
\proof
The statement follows from Koszul formula making usual computations.

\gata

Having determined Levi Civita connection, we can compute now the Riemannian
curvature tensor $\tilde{R}^a$ on $T(M)$. We give

\begin{prop}
The curvature tensor is given by

$\begin{array}{l}
\quad\tilde R^a_{X^HY^H}Z^H=(R_{XY}Z)^H+\frac a4\ \left[
   R_{\uu R_{XZ}\uu}Y-R_{\uu R_{YZ}\uu}X+2R_{\uu R_{XY}\uu}Z
    \right]^H+\\
    \qquad\qquad\qquad +\frac 12\ \left[
   (\nabla_ZR)_{XY}\uu
    \right]^V
\end{array}
    $\\[2mm]
$\begin{array}{l}
\quad\tilde R^a_{X^HY^H}Z^V=\left[ R_{XY}Z+
  \frac a4\ (R_{YR_{\uu Z}X}\uu-R_{XR_{\uu Z}Y}\uu)
    \right]^V+\LL g(Z,\uu)(R_{XY}\uu)^V+\\
   \qquad\qquad\qquad +\frac{1-\LL}{\rr^2}\ g(R_{XY}\uu,Z)\uu^V+\frac a2\ \left[
(\nabla_XR)_{\uu Z}Y-(\nabla_YR)_{\uu Z}X
    \right]^H
    \end{array}
    $\\[2mm]
$\begin{array}{l}
\quad \tilde R^a_{X^HY^V}Z^H=\frac a2\left[(\nabla_XR)_{\uu Y}Z\right]^H+\\
\qquad\qquad
+\frac12\left[
    R_{XZ}Y-\frac a2\ R_{XR_{\uu Y}Z}\uu+\LL g(Y,\uu)R_{XZ}\uu+\frac{1-\LL}{\rr^2}\ g(R_{XZ}\uu,Y)\uu
    \right]^V
 \end{array}
    $\\[2mm]
$\begin{array}{l}
\quad\tilde R^a_{X^HY^V}Z^V=-\frac a2\ (R_{YZ}X)^H-\frac{a^2}4\ (R_{\uu Y}R_{\uu Z}X)^H+\\
\qquad\qquad\qquad
 +\frac{a'}4\left[ g(Z,\uu)(R_{\uu Y}X)^H-g(Y,\uu)(R_{\uu Z}X)^H
 \right]
 \end{array}
    $
\begin{equation}
\label{eq:tildeRA}
\begin{array}{l}
\tilde R^a_{X^VY^V}Z^H=a (R_{XY}Z)^H+\frac{a'}2\left[
  g(X,\uu)R_{\uu Y}Z-g(Y,\uu)R_{\uu X}Z
    \right]^H\qquad\quad\ \\
    \qquad\qquad\qquad+\frac {a^2}4\left[
 R_{\uu X}R_{\uu Y}Z-R_{\uu Y}R_{\uu X}Z
\right]^H
\end{array}
\end{equation}
$\begin{array}{l}
\quad\tilde R^a_{X^VY^V}Z^V=F_1(t) g(Z,\uu)\left[g(X,\uu)Y^V-g(Y,\uu)X^V\right]+\\
      \qquad\qquad\qquad +F_2(t)\left[g(X,Z)Y^V-g(Y,Z)X^V\right]+\\
      \qquad\qquad\qquad + F_3(t)\left[g(X,Z)g(Y,\uu)-g(Y,Z)g(X,\uu)\right]\uu^V,
\end{array}
  $\\[2mm]
where $F_1=\LL'+\frac{\LL(1-\LL)}{\rr^2}$, $F_2=\LL^2-\frac{(1-\LL)^2}{\rr^2}$ and
$F_3=\frac{\LL'-\LL^2}{\rr^2}+\frac{1-\LL}{\rr^4}\ $.
\end{prop}
In the following let $\tilde Q^a(U,V)$ denote the square of the area of the parallelogram with
sides $U$ and $V$ for $U,V\in\chi(T(M))$,
$$\tilde Q^a(U,V)=g_a(U,U)g_a(V,V)-g_a(U,V)^2.
$$
We have
\begin{lemma}
Let $X,Y\in T_pM$ be two orthonormal vectors. Then
\begin{equation}
\left\{
\begin{array}{l}
\tilde Q^a(X^H,Y^H)=1\ ,\
\tilde Q^a(X^H,Y^V)=a(t)\big(1+g(Y,\uu)^2\big)\\[2mm]
\tilde Q^a(X^V,Y^V)=a(t)^2\big(1+g(X,\uu)^2+g(Y,\uu)^2\big).
\end{array}
\right.
\end{equation}
\end{lemma}
We compute now the sectional curvature of the Riemannian manifold $(T(M), g_a)$, namely
$\tilde K^a(U,V)=\frac{g_a(\tilde R^a_{UV}V,U)}{\tilde Q^a(U,V)}\ $
for $U,V\in\chi(T(M))$.\\[2mm]
Denote by $T_0(M)=T(M)\setminus{\mathbf{0}}$ the tangent bundle of non-zero vectors tangent to $M$.\\[2mm]
For a given point $(p,\uu)\in T_0(M)$ consider an orthonormal basis $\{e_i\}_{i=\overline{1,m}}$ for the tangent space
$T_p(M)$ of $M$ such that $e_1=\frac u{|u|}$. Consider on $T_{(p,\uu)}T(M)$ the following vectors
\begin{equation}
\label{eq:E}
\begin{array}{l}
E_i=e_i^H,\ i=\overline{1,m} \ , \
E_{m+1}=\frac1{\rr\sqrt{a}}\ e_1^V\ ,\
E_{m+k}=\frac 1{\sqrt{a}}\ e_k^V,\ k=\overline{2,m}.
\end{array}
\end{equation}
It is easy to check that $\{E_1,\ldots,E_{2m}\}$ is an orthonormal basis in $T_{(p,\uu)}T(M)$ (with respect to the metric $g_a$).
We will write the expressions of the sectional curvature $\tilde K^a$ in terms of this basis. We have
\begin{equation}
\left\{
\begin{array}{l}
\tilde K^a(E_i,E_j)=K(e_i,e_j)-\frac {3a(t)}4\ \left|R_{e_ie_j}\uu\right|^2,\\
\tilde K^a(E_i,E_{m+1})=0, \\
\tilde K^a(E_i,E_{m+k})=\frac 14\left|R_{\uu e_k}e_i\right|^2, \\
\tilde K^a(E_{m+1}E_{m+k})=-\frac{F_2+2tF_3}{a(t)}, \\
\tilde K^a(E_{m+k}E_{m+l})=-\frac{F_2}{a(t)},\quad i,j=1,\ldots,m;\ k,l=2,\ldots,m.
\end{array}
\right.
\end{equation}
Here $|\cdot|$ denotes the norm of the vector with respect to the metric $g$ (in a point).\\[2mm]
{\bf Question:} {\em Can we have constant sectional curvature $\tilde c$ on $T(M)$?}\\[2mm]
If this happens, then it must be $0$, so $T(M)$ is flat.
One gets easily that $M$ is locally Euclidean.
Then, we should also have $F_2(t)=0$. It follows, $F_3(t)=0$ and $F_1(t)=0$.
On the other hand an ordinary differential equation occurs:\\[2mm]
\centerline{$\frac{a'(t)}{a(t)}=\frac2{1+\sqrt{1+2t}}$.}\\[2mm]
A simple computation shows that
\begin{equation}
a(t)=a_0 \frac{e^{2\sqrt{1+2t}}}{(1+\sqrt{1+2t})^2}\quad ,\  a_0>0.
\end{equation}
\begin{remark}\rm
The manifold $T(M)$ equipped with the Cheeger Gromoll has non constant sectional curvature.
\end{remark}

Putting $a_0$, such that $a(0)=1$ we can state the following
\begin{teo}
\label{prop:2.17}
Consider
$g_1$ on $T(M)$ given by
\begin{equation}
\label{metric:g1}
\left\{\begin{array}{l}
g_1(X^H,Y^H)=g(X,Y),\ g_1(X^H,Y^V)=0\\[2mm]
g_1(X^V,Y^V)=\frac{4e^{2(\rr-1)}}{(1+\rr)^2}\ \big(g(X,Y)+g(X,\uu)g(Y,\uu)\big)
\end{array}\right.
\end{equation}

The manifold $(T(M),g_1)$ is flat.
\end{teo}
Let us now compare the scalar curvatures of $(M,g)$ and $(T(M),g_a)$.
\begin{prop}
Let $(M,g)$ be a Riemannian manifold and endow the tangent bundle $T(M)$ with the metric $g_a$. Let
${\rm scal}$ and $\widetilde{\rm scal}^a$ be the scalar curvatures of $g$ and $g_a$ respectively.
The following relation holds:
\begin{equation}
\widetilde{\rm scal}^a={\rm scal}+\frac{2-3a}2\ \sum\limits_{i<j}\left|R_{e_ie_j}\uu\right|^2+
    \frac{1-m}a\ \left(m F_2+4t F_3\right),
\end{equation}
where $\{e_i\}_{i=1,\ldots,m}$ is a local orthonormal frame on $T(M)$.
\end{prop}
\proof
Using that ${\rm scal}=\sum\limits_{i\neq j} K(e_i,e_j)$ and the formula\\[1mm]
\centerline{
$\sum\limits_{i,j=1}^m\left|R_{e_i\uu}e_j\right|^2=\sum\limits_{i,j=1}^m\left|R_{e_ie_j}\uu\right|^2$
}\\[1mm]
we get the conclusion.
\gata\\[2mm]
Consider $M$ a real space form with $c$ the constant sectional curvature.\\[2mm]
{\bf Question:}
{\em Could we find functions $a$ such that $T(M)$ equipped with the metric $g_a$
has constant scalar curvature?}\\[2mm]
Then $(T(M),g_a)$
has constant scalar curvature if and only if $a$ satisfies the following ODE:
$$
\begin{array}{l}
-\frac{1}{2\  {{(1+2\  t)}^2}\  {{a(t)}^3}}  \\
\noalign{\vspace{0.5mm}}
\hspace{2.em} \big(-2\  (m+2\  (-2+m)\  t)\  {{a(t)}^2}-4\  t\  {{(c+2\  c\  t)}^2}\  {{a(t)}^3}+  \\
\noalign{\vspace{0.5mm}}
\hspace{4.em} +6\  t\  {{(c+2\  c\  t)}^2}\  {{a(t)}^4}+(-6+m)\  t\  (1+2\  t)\  {{{a^{\prime }}(t)}^2}+
 \\
\noalign{\vspace{0.5mm}}
\hspace{4.em} +2\  a(t)\  ((m+2\  (-1+m)\  t)\  {a^{\prime }}(t)+2\  t\  (1+2\  t)\  {a^\prime{}^\prime}(t))\big)
    ={\rm const.}
\end{array}
$$
which seems to be very complicated to solve.\\[3mm]
{\bf Acknowledgements.} This work was partially supported by
Grant CNCSIS 1463/ n.18/2005.

\smallskip
\begin{minipage}{2.5in}
\begin{flushleft}
Marian Ioan MUNTEANU\\
Faculty of Mathematics \\
Al.I.Cuza University of Ia\c si,\\
Bd. Carol I, n. 11\\
700506 - Ia\c si, ROMANIA \\
e-mail: munteanu@uaic.ro
\end{flushleft}
\end{minipage}

\end{document}